\documentstyle[twoside]{article}
\def\C{{\ifmmode{C\hskip-5.0pt
                  \vrule height6.0pt depth 0pt width0.6pt\hskip5pt}
   \else{\hbox{$C\hskip-5.0pt
                  \vrule height6.0pt depth 0pt width0.6pt\hskip5pt$}}\fi}}

\def\beq{\begin{equation}}
\def\eeq{\end{equation}}
\def\bea{\begin{eqnarray}}
\def\eea{\end{eqnarray}}

\def\ba{\begin{array}}
\def\ea{\end{array}}

\def\titem#1{\par\hang\indent\hbox to 0pt{\hss\hbox to \parindent{
\enspace\hss#1}}\ignorespaces}

\def\mod{\mathop{\rm mod}}
\def\sds{\mathbin{\rlap{\raise0.9pt\hbox{$\scriptstyle+$}}\!\!\supset}}
\def\anystrut#1{\dimen11=#1\dimen12=#1
\divide\dimen12 by 4\dimen11=\dimen12\multiply\dimen11 by 3
\hbox{\vrule height\dimen11 depth\dimen12 width0pt}}
\newdimen\tadhdimen
\newdimen\tabhdimen
\edef\tabstrut{\anystrut{12truept}} 

\def\youngt#1{
\vcenter{\offinterlineskip
\halign{&\tabstrut\hbox to \tabhdimen{\hss$##$\hss}\cr #1}}}
\def\youngd#1{
\vcenter{\offinterlineskip
\halign{&\vrule##&\tabstrut\hbox to \tadhdimen{\hss$##$\hss}\cr #1}}}

\newdimen\stadhdimen
\newdimen\stabhdimen
\edef\shortstrut{\anystrut{6.6truept}} 

\def\syoungd#1{{\vcenter{\offinterlineskip
\halign{&\vrule##&\shortstrut
\hbox to \stabhdimen{\hss$\scriptstyle ##$\hss}\cr #1}}}}
\def\syoungt#1{{\vcenter{\offinterlineskip
\halign{&\shortstrut
\hbox to \stadhdimen{\hss$\scriptstyle ##$\hss}\cr #1}}}}

\begin{document}
\renewcommand{\thefootnote}{\fnsymbol{footnote}}
\setcounter{footnote}{1}

\begin{center}
{\bf TWO-ROWED HECKE ALGEBRA REPRESENTATIONS AT ROOTS OF UNITY}%
\footnote{Presented at the 4th
International Colloquium ``Quantum Groups and Integrable Systems,''
Prague, 22-24 June 1995; and to appear in proceedings in
{\it Czech J.~Phys.}}

\vspace{1cm}

{\sc Trevor Alan Welsh}\footnote{
E-mail~: taw@maths.soton.ac.uk}\\[2mm]
{\sl Faculty of Mathematical Studies,} \\
{\sl University of Southampton, Southampton, SO17 1BJ, U.K.}
\end{center}

\vspace{1cm}

{\small
In this paper, we initiate a study into the explicit construction of
irreducible representations of the Hecke algebra $H_n(q)$ of type
$A_{n-1}$ in the non-generic case where $q$ is a root of unity.
The approach is via the Specht modules of $H_n(q)$ which are
irreducible in the generic case, and possess a natural basis indexed
by Young tableaux.
The general framework in which the irreducible non-generic
$H_n(q)$-modules are to be constructed is set up and, in particular,
the full set of modules corresponding to two-part partitions is
described. Plentiful examples are given.}

\vspace{5mm}
\begin{center}
{\bf 1\ \ Introduction and notation}
\end{center}

The Hecke algebra $H_n(q)$ (of type $A_{n-1}$) is the unital associative
algebra over $\C$, generated by $h_i$, $i=1,2,\ldots,n-1$, subject
to the relations:
\beq
\ba{l}
h_ih_{i+1}h_i=h_{i+1}h_ih_{i+1};\\
h_ih_j=h_jh_i\qquad\hbox{for }\vert i-j\vert>1;\\
h_i^2=(q-1)h_i+q.
\label{hoopy}
\ea\eeq

The parameter $q\in\C$ will be permitted to take any non-zero value.
It is said to be generic if $q=1$ or $q^p\ne 1$ for $p=2,3,4,\ldots$.
Otherwise, if $q$ is a primitive $p$th root of unity for $p>2$,
it is said to be non-generic.

When $q=1$, $H_n(q)$ may be identified with the group algebra
$\C S_n$ of the symmetric group on $n$ symbols, through identifying
each $h_i$ with the simple transposition $s_i=(i,i+1)\in S_n$.

If $w=s_{i_1}s_{i_2}\cdots s_{i_k}$ and $w\in S_n$ cannot be
expressed as a shorter product of the generators $s_i$, then
$s_{i_1}s_{i_2}\cdots s_{i_k}$ is said to be a reduced expression
for $w$ and the value of $k$ is the length $l(w)$ of $w$.
Thereupon, the relations (\ref{hoopy}) imply that the map
$h:\C S_n\rightarrow H_n(q)$ for which $h(s_i)=h_i$ and
$h(ww^\prime)=h(w)h(w^\prime)$ for $w,w^\prime\in S_n$ satisfying
$l(ww^\prime)=l(w)+l(w^\prime)$, and extended linearly, is well defined.
It follows that if $l(w)=k$ and $w=s_{i_1}s_{i_2}\cdots s_{i_k}$,
then $h(w)=h_{i_1}h_{i_2}\cdots h_{i_k}$. Furthermore, a basis of
$H_n(q)$ is provided by $\{h(w):w\in S_n\}$.
\par
It may be shown that if $q$ is generic then $H_n(q)$
is isomorphic to $\C S_n$ \cite{DJ86,Wn88} and
the representation theory of $H_n(q)$ is much the
same as that of $S_n$. In particular, the inequivalent irreducible
representations of $H_n(q)$ are indexed by partitions $\lambda$ of $n$.
That is, by finite integer sequences
$\lambda=(\lambda_1,\lambda_2,\ldots,\lambda_r)$ for which
$\lambda_1+\lambda_2+\cdots+\lambda_r=n$ and
$\lambda_1\ge\lambda_2\ge\cdots\ge\lambda_r>0$.
A partition for which no part $\lambda_i$
is repeated more than $p-1$ times is said to be $p$-regular.
In Section 2, an explicit construction of
the irreducible modules of $H_n(q)$ with $q$ generic will be described.
This generalisation of the well known Specht module
construction (see \cite{JK81}) was first described in \cite{KWy92},
and is based on the use of Young diagrams, Young tableaux
and $q$-analogues of Young symmetrisers.
The Young diagram $F^\lambda$ associated with the partition
$\lambda=(\lambda_1,\lambda_2,\ldots,\lambda_r)$ is a left-adjusted,
top-adjusted array of square boxes such that the $i$th row
(counting from the top) contains $\lambda_i$ boxes.
For instance, if $\lambda=(5,3,2,2)$, then
\beq
F^\lambda=\syoungd{
\multispan{11}\hrulefill\cr
&&&&&&&&&&\cr \multispan{11}\hrulefill\cr
&&&&&&\cr \multispan{7}\hrulefill\cr
&&&&\cr \multispan{5}\hrulefill\cr
&&&&\cr \multispan{5}\hrulefill\cr}.
\eeq
Filling (or replacing) each of the $n$ boxes of $F^\lambda$ with distinct
elements of $\{1,2,\ldots,n\}$ yields what is known as a Young tableau.
Of the possible $n!$ tableaux of a given shape,
those for which the entries are increasing across each row and down
each column are
known as {\em standard\/} tableaux.
Examples may be found at (\ref{std1}), (\ref{std2}) and (\ref{std3}).
That particular standard tableau of shape $\lambda$ for which the entries
increase down first the leftmost column and then down successive columns
taken left to right (e.g.~(\ref{std1})) is denoted $t^\lambda_-$.
The number of standard tableaux of shape $\lambda$ is equal to
the dimension of the irreducible representation of $S_n$
(and $H_n(q)$ with $q$ generic) labelled by $\lambda$ (see \cite{JK81}).
In fact, the Specht module construction enables a basis to be
identified naturally with the set of standard tableaux.

\begin{center}
{\bf 2\ \ The Specht modules}
\end{center}

If $\lambda$ is a partition of $n$, the Specht module
$S^\lambda$ of $H_n(q)$ is defined to be the linear span of the vectors
$v_{t^\lambda}$, indexed by Young tableaux $t^\lambda$
and subject to certain relations (which will be defined below).
The {\em natural\/} action of $H_n(q)$ on these vectors is defined in
the following way.
We say that the entry $i$ precedes $j$ in $t^\lambda$ if $i$ occurs
before $j$ on reading the entries of $t^\lambda$
down the first and then successive columns.
If $x^\lambda$ is identical to $t^\lambda$ apart from the transposition
of $i$ and $i+1$, then $h_i$ acts on $v_{t^\lambda}$ as follows:
\beq
h_i v_{t^\lambda} =
\cases{
v_{x^\lambda}
&if $i$ precedes $i+1$ in $t^\lambda$;\cr
q v_{x^\lambda} + (q-1) v_{t^\lambda}
&if $i+1$ precedes $i$ in $t^\lambda$.\cr}
\label{action}
\eeq
\par
It is possible to express every $v_{z^\lambda}$ in terms of standard
tableaux, by means of the following two types of relation:
\smallskip\begingroup\parindent=22pt
\titem{1. } {\bf Column relations.} Entries within a column may be
transposed, if the corresponding vector is multiplied by $-1$.
Thus if $x^\lambda$ differs from $z^\lambda$ only in that a single pair of
entries within a column are transposed then:
\beq
v_{z^\lambda}=-v_{x^\lambda}.
\label{colrel}
\eeq
For example (denoting $v_{t^\lambda}$ by $t^\lambda$ for
typographical reasons),
\beq
\youngt{1&8&5&10&4&12\cr 6&11&3\cr 9&2&7\cr 13\cr}
=-\youngt{1&8&5&10&4&12\cr 6&2&3\cr 9&11&7\cr 13\cr}
=-\youngt{1&2&3&10&4&12\cr 6&8&5\cr 9&11&7\cr 13\cr}.
\eeq
\titem{2. } {\bf Garnir relations.} Assume that $z^\lambda$ is such that
its entries increase down each column.
If $z^\lambda$ is not standard then an adjacent pair of entries
exists with that on
the left greater than that on the right. Consider these two entries
together with all those below the left one and all those above the right one.
For example, we could consider the highlighted entries in:
\beq
\youngt{1&2&\bf3&10&4&12\cr 6&\bf8&\bf5\cr 9&\bf11&7\cr 13\cr}.
\eeq
Now form all possible tableaux $t^\lambda$ by permuting these entries
in all ways such that the permuted entries are increasing
down the portions of each of the two columns being considered.
The Garnir relation is then the following expression in which the sum
is over all such tableaux:
\beq
(-q)^{l(w_{z^\lambda})}
\sum_{t^\lambda} (-q)^{-l(w_{t^\lambda})}
v_{t^\lambda}=0,
\label{garnirrel}
\eeq
where $w_{t^\lambda}\in S_n$ maps $t^\lambda_-$ to $t^\lambda$.
The above example gives the
Garnir relation:
\beq
\ba{l}
\!\!\youngt{1&2&\bf3&10&4&12\cr 6&\bf8&\bf5\cr 9&\bf11&7\cr 13\cr}
-q\youngt{1&2&\bf3&10&4&12\cr 6&\bf5&\bf8\cr 9&\bf11&7\cr 13\cr}
+q^2\youngt{1&2&\bf5&10&4&12\cr 6&\bf3&\bf8\cr 9&\bf11&7\cr 13\cr}\\
+q^2\youngt{1&2&\bf3&10&4&12\cr 6&\bf5&\bf11\cr 9&\bf8&7\cr 13\cr}
-q^3\youngt{1&2&\bf5&10&4&12\cr 6&\bf3&\bf11\cr 9&\bf8&7\cr 13\cr}
+q^4\youngt{1&2&\bf8&10&4&12\cr 6&\bf3&\bf11\cr 9&\bf5&7\cr 13\cr}=0.\\
\ea
\label{garnireg}
\eeq
\endgroup
\par\noindent
As in the example above, these relations do not necessarily immediately
express an arbitrary $v_{t^\lambda}$ in terms of standard tableaux.
However, it may be shown through employing a suitable order on the set
of tableaux \cite{JK81}, that repeated application of the column
and Garnir relations enables any term to be rendered in terms of
standard tableaux in a finite number of steps.
This completes the construction of the irreducible Specht module
$S^\lambda$ of $H_n(q)$ since the number of standard tableaux is equal
to the dimension of the representation of $H_n(q)$ indexed by $\lambda$
and consequently,
\beq
\{v_{t^\lambda}:t^\lambda {\rm\ is\ standard}\}
\eeq
is a basis for $S^\lambda$.

As an example, consider representing $h_1\in H_5(q)$ in the Specht
module $S^{(3,2)}$, by acting on each $v_{t^{(3,2)}}$ for which $t^{(3,2)}$
is standard (once more $v_{t^\lambda}$ is written as $t^\lambda$):
$$
\tabhdimen=8.22truept
\ba{l}
{ h_1} \youngt{1&3&5\cr 2&4\cr}
= \youngt{2&3&5\cr 1&4\cr}
= -\youngt{1&3&5\cr 2&4\cr},\\
{ h_1} \youngt{1&2&5\cr 3&4\cr}
= \youngt{2&1&5\cr 3&4\cr}
= q\youngt{1&2&5\cr 3&4\cr}-q^2\youngt{1&3&5\cr 2&4\cr},\\
{ h_1} \youngt{1&3&4\cr 2&5\cr}
= \youngt{2&3&4\cr 1&5\cr}
= -\youngt{1&3&4\cr 2&5\cr},\\
{ h_1} \youngt{1&2&4\cr 3&5\cr}
= \youngt{2&1&4\cr 3&5\cr}
= q\youngt{1&2&4\cr 3&5\cr}-q^2\youngt{1&3&4\cr 2&5\cr},\\
{ h_1} \youngt{1&2&3\cr 4&5\cr}
= \youngt{2&1&3\cr 4&5\cr}
= q\youngt{1&2&3\cr 4&5\cr}-q^2\youngt{1&4&3\cr 2&5\cr}
=q \youngt{1&2&3\cr 4&5\cr}
-q^3 \youngt{1&3&4\cr 2&5\cr}
+q^4 \youngt{1&3&5\cr 2&4\cr}.\\
\ea
$$
Here, column relations have been used in the first and third calculations,
and Garnir relations have been used in the second, fourth and last (twice),
to express each result in terms of the standard tableaux.
Consequently, in the representation labelled by the partition $(3,2)$,
$h_1$ is mapped to the matrix (where zeros are denoted by dots):
\beq
\pmatrix{
-1&-q^2&.&.&q^4\cr
.&q&.&.&.\cr
.&.&-1&-q^2&-q^3\cr
.&.&.&q&.\cr
.&.&.&.&q\cr}.
\eeq
The matrices representing the generators $h_i$ of $H_n(q)$ in each
irreducible representation for $n\le5$ given
in \cite{KWy92} have been produced in a similar way.

\begin{center}
{\bf 3\ \ The Young symmetriser and its annihilators}
\end{center}

For each entry $a$ of $t^\lambda_-$ which is not at the bottom of
a column, define the {\em column element}:
\beq
C^\lambda_a=1+h_a.
\label{coldef}
\eeq
Its action on $v_{t^\lambda_-}$ gives rise to a Column relation
(cf.~(\ref{colrel})):
\beq
C^\lambda_a v_{t^\lambda_-}=0.
\label{annih1}
\eeq

The {\em Garnir element\/} $G^\lambda_a$ is defined for each $a$ which is
not at the end of a row of $t^\lambda_-$, through first letting $d$ be the
entry to the right of $a$, $b$ be the entry at the bottom of the column
containing $a$, and $c$ (${}=b+1$) the entry at the top of the column
containing $d$ in $t^\lambda_-$.
With $W_{ij}$ the subgroup of $S_n$ permuting only $\{i,i+1,\ldots,j\}$,
let ${\cal G}^\lambda_a$ be a set of
left coset representatives for $W_{ab}\times W_{cd}$ in $W_{ad}$
chosen so that each representative is of minimal length in its coset
(it is unique). Then let \cite{KWy92}:
\beq
G^\lambda_a=\sum_{d\in{\cal G}^\lambda_a} (-q)^{-l(d)} h(d).
\label{garnirdef}
\eeq
Its action on $v_{t^\lambda_-}$ gives rise to a Garnir relation
(cf.~(\ref{garnirrel})):
\beq
G^\lambda_a v_{t^\lambda_-}=0,
\label{annih2}
\eeq

It is easily shown that the general column and Garnir relations
of Section 2 are a consequence of (\ref{annih1}) and (\ref{annih2}).
These properties themselves arise by identifying
$v_{t^\lambda_{\smash{-}}}$
with the $q$-analogue $Y^\lambda(q)$ of the Young symmetriser.
$Y^\lambda(q)$ was originally defined in \cite{DJ86,Gy86} and cast
in a form suitable for the current purposes in \cite{KWy92,BKW95}.
However, as is seen, only its $2n-r-\lambda_1$ annihilators $C_a$ and $G_a$
are required in the construction of the Specht module $S^\lambda$.
Thus $S^\lambda$ may be defined as the free module generated by a non-zero
vector
(say $v_{t^\lambda_{\smash{-}}}$)
subject to (\ref{annih1}) and (\ref{annih2}).
This viewpoint of $S^\lambda$ will be utilised in what follows.
\par
To illustrate it, consider $\lambda=(6,3,3,1)$, for which:
\beq
t^\lambda_-=
\youngt{1&5&8&11&12&13\cr 2&6&9\cr 3&7&10\cr 4\cr}.
\label{std1}
\eeq
Here we have the seven column elements $1+h_1$, $1+h_2$, $1+h_3$,
$1+h_5$, $1+h_6$, $1+h_8$ and $1+h_9$, each of which annihilates
$v_{t^\lambda}$. There are nine Garnir elements $G^\lambda_a$ for
$a=1,2,3,5,6,7,8,11,12$, each of which annihilates $v_{t^\lambda}$.
Typically:
\beq
\ba{l}
G^\lambda_6=
q^4-q^3h_7+q^2h_6h_7+q^2h_8h_7-qh_6h_8h_7+h_7h_6h_8h_7;\\
G^\lambda_8=
q^3-q^2h_{10}+qh_9h_{10}-h_8h_9h_{10};\\
G^\lambda_{11}=
q-h_{11}.\\
\ea\eeq
In fact $G^\lambda_6 v_{t^\lambda_-}=0$ gives rise to (\ref{garnireg}).

\begin{center}
{\bf 4\ \ Decomposing \boldmath $S^\lambda$ at roots of unity}
\end{center}

In the generic case when $q$ is not a root of unity, each Specht module
$S^\lambda$ of $H_n(q)$ is irreducible.
However, this is no longer so if
$q$ is a root of unity, although $S^\lambda$ remains well-defined.
For $q$ a primitive $p$th root of unity,
let $D^\lambda_p$ be the irreducible $H_n(q)$-module
obtained by factoring out the maximal proper submodule from $S^\lambda$.
It is shown in \cite{DJ86} that in this case,
\beq
\{D^\lambda_p: \lambda {\rm\ is\ } p {\rm -regular}\}
\eeq
is a complete set of irreducible, irredundant $H_n(q)$-modules.
Very little is known about the $D^\lambda_p$ or the composition
series of $S^\lambda$ in terms of the $D^\lambda_p$ except in a few
specific cases (see \cite{Jm90} for $n\le10$ and \cite{CK92} for $n\le5$).
\par
The viewpoint developed in the previous section provides a means of
tackling these questions in a quite general way.
It relies on the fact that $D^\mu_p$ is characterised by the
presence of a non-zero vector
$v_{t^\mu_{\smash{-}}}$
which is annihilated by the set of column and Garnir elements,
$C^\mu_a$ and $G^\mu_a$ defined above.
This follows because, via (\ref{action}),
$v_{t^\mu_{\smash{-}}}$
generates the whole of $S^\mu$, and hence
$v_{t^\mu_{\smash{-}}}$
cannot be present in any proper submodule.
Thus, to determine whether $S^\lambda$ is reducible, it is sufficient
to prove the existence of a non-zero $v^\mu\in S^\lambda$ having
the same set of annihilators as
$v_{t^\mu_{\smash{-}}}\in S^\mu$
for some $p$-regular partition $\mu\ne\lambda$ of $n$.
Conversely, the absence of all such
$v_{t^\mu_{\smash{-}}}$
would prove $S^\lambda$ to be irreducible.
(In fact, the results of \cite{DJ86} and \cite{DJ87} considerably restrict
the set of $\mu$ for which $D^\mu_p$ may occur as a composition factor
of $S^\lambda$.)
\par
\begingroup\tabhdimen=8.22truept
As an example, consider $\lambda=(3,2)$.
We will show that if $p=3$ then
\beq
v^\mu=(1+h_4)v_{t^{(3,2)}_-}=
\youngt{1&3&5\cr 2&4\cr}+\youngt{1&3&4\cr 2&5\cr}
\label{anniheg1}
\eeq
is annihilated by the column and Garnir elements of $\mu=(4,1)$, and hence
that $S^{(3,2)}$ has a submodule $D^{(4,1)}_3$.
Since
$t^{(4,1)}_-=\syoungt{1&3&4&5\cr 2\cr}$,
the column and Garnir elements of $\mu=(4,1)$, are:
\beq
\ba{ll}
\romannumeral1) & (1+h_1); \\
\romannumeral2) & (q^2-qh_2+h_1h_2); \\
\romannumeral3) & (q-h_3); \\
\romannumeral4) & (q-h_4). \\
\ea\eeq
Acting on (\ref{anniheg1}) with each of these, using (\ref{action}) gives:
$$
\ba{ll}
\romannumeral1) & (1+h_1)v^\mu=
\youngt{1&3&5\cr 2&4\cr}+\youngt{2&3&5\cr 1&4\cr}
+\youngt{1&3&4\cr 2&5\cr}+\youngt{2&3&4\cr 1&5\cr}
=0;\\
\romannumeral2) & (q^2-qh_2+h_1h_2) v^\mu=
q^2\youngt{1&3&5\cr 2&4\cr}-q\youngt{1&2&5\cr 3&4\cr}+\youngt{2&1&5\cr 3&4\cr}
\\
&\qquad\qquad\qquad\qquad\qquad\qquad\qquad\qquad\qquad
+q^2\youngt{1&3&4\cr 2&5\cr}-q\youngt{1&2&4\cr 3&5\cr}+\youngt{2&1&4\cr 3&5\cr}
=0;\\
\romannumeral4) & (q-h_4)v^\mu=
q\youngt{1&3&5\cr 2&4\cr}-\youngt{1&3&4\cr 2&5\cr}
+q\youngt{1&3&4\cr 2&5\cr}-q\youngt{1&3&5\cr 2&4\cr}
-(q-1)\youngt{1&3&4\cr 2&5\cr}
=0;\\
\romannumeral3) & (q-h_3)v^\mu=
q\youngt{1&3&5\cr 2&4\cr}-\youngt{1&4&5\cr 2&3\cr}
+q\youngt{1&3&4\cr 2&5\cr}-\youngt{1&4&3\cr 2&5\cr}
\\
&\qquad
=(1+q)\youngt{1&3&5\cr 2&4\cr}+q\youngt{1&3&4\cr 2&5\cr}
-q\youngt{1&3&4\cr 2&5\cr}+q^2\youngt{1&3&5\cr 2&4\cr}
=(1+q+q^2)\youngt{1&3&5\cr 2&4\cr}=0,
\\
\ea
$$
since if $p=3$, then $1+q+q^2=0$.
Therefore, $D^{(4,1)}_3$ is a submodule of $S^{(3,2)}$.
It may be shown that the 4-dimensional $S^{(4,1)}$ is irreducible when
$p=3$, so that $D^{(4,1)}_3\equiv S^{(4,1)}$.
Hence $D^{(3,2)}_3$ is of dimension $5-4=1$.
It is spanned by $v_{t^{(3,2)}_-}$.
\endgroup

In order to express the general structure in the case of two-part
partitions, we introduce the notion of a boundary strip
of a Young diagram $F^\lambda$.
It is a continuous strip of boxes obtained by starting at the
rightmost end of a row of $F^\lambda$ and, for a number of steps,
recursively passing to the box below if one exists, otherwise passing
to the box to the left. It must end at the bottom of a column.
The length of the boundary strip is the number of boxes it comprises.

\proclaim Theorem 1. If $\lambda=(\lambda_1,\lambda_2)$ and $q$
is a primitive $p$th root of unity then
$S^\lambda$ is reducible if and only if for some integer $k>0$,
$F^\lambda$ has a boundary strip of length $kp$ having at least one,
but not more than $p-1$ boxes in the second row
(or equivalently, if there exists an integer $k>0$ such that
$\lambda_1-\lambda_2+2\le kp\le\min\{\lambda_1+1,\lambda_1-\lambda_2+p\}$).
If so, $S^\lambda$ has an irreducible submodule corresponding to the diagram
obtained by moving all the boxes of the boundary strip into the top row.
The corresponding quotient module is irreducible.
That is:
\beq
S^{(\lambda_1,\lambda_2)}=
D^{(\mu_1,\mu_2)}_p
\sds D^{(\lambda_1,\lambda_2)}_p,
\label{reduction}
\eeq
where
\beq
(\mu_1,\mu_2)=(\lambda_1+p-1-(\lambda_1-\lambda_2)\mod p,
\lambda_2-p+1+(\lambda_1-\lambda_2)\mod p)
\eeq
\par\noindent
This theorem is illustrated by the following table, which for various
$\lambda=(\lambda_1,\lambda_2)$ and $p$, displays the Young diagram $F^\lambda$
with the appropriate boundary strip indicated, says whether $S^\lambda$
is reducible, and shows its composition series.
$$
\vcenter{\openup\jot
\halign{\enspace$#$\quad\hfill&$#$\qquad\hfill&$#$\qquad\hfill&#\quad\hfill
&$#$\enspace\hfill\cr
\enspace\lambda&p&F^\lambda&$S^\lambda$&\rm{Composition}\cr
(5,4)&3&
\syoungd{\multispan{11}\hrulefill\cr
 &&&&&&&\bullet&&\bullet&\cr \multispan{11}\hrulefill\cr
 &&&&&&&\bullet&\cr \multispan{9}\hrulefill\cr}
&reducible&S^{(5,4)}=D^{(6,3)}_3\sds D^{(5,4)}_3\cr
(6,4)&3&
\syoungd{\multispan{13}\hrulefill\cr
 &&&&&&&&&&&&\cr \multispan{13}\hrulefill\cr
 &&&&&&&&\cr \multispan{9}\hrulefill\cr}
&irreducible&S^{(6,4)}=D^{(6,4)}_3\cr
(7,4)&3&
\syoungd{\multispan{15}\hrulefill\cr
 &&&&&&&\bullet&&\bullet&&\bullet&&\bullet&\cr \multispan{15}\hrulefill\cr
 &&&&&\bullet&&\bullet&\cr \multispan{9}\hrulefill\cr}
&reducible&S^{(7,4)}=D^{(9,2)}_3\sds D^{(7,4)}_3\cr
(8,3)&3&
\syoungd{\multispan{17}\hrulefill\cr
 &&&&&&&&&&&&&&&&\cr \multispan{17}\hrulefill\cr
 &&&&&&\cr \multispan{7}\hrulefill\cr}
&irreducible&S^{(8,3)}=D^{(8,3)}_3\cr
(8,2)&5&
\syoungd{\multispan{17}\hrulefill\cr
 &&&&&&&&&&&&&&&&\cr \multispan{17}\hrulefill\cr
 &&&&\cr \multispan{5}\hrulefill\cr}
&irreducible&S^{(8,2)}=D^{(8,2)}_5\cr
(9,3)&5&
\syoungd{\multispan{19}\hrulefill\cr
 &&&&&\bullet&&\bullet&&\bullet&&\bullet&&\bullet&&\bullet&&\bullet&\cr
 \multispan{19}\hrulefill\cr
 &\bullet&&\bullet&&\bullet&\cr \multispan{7}\hrulefill\cr}
&reducible&S^{(9,3)}=D^{(12)}_5\sds D^{(9,3)}_5\cr }}
$$
\par
Theorem 1 has the consequence that the character $\tilde\chi^\lambda_p$
of the irreducible representation corresponding to $D^\lambda_p$
may be expressed as a finite sum over the characters
$\chi^\lambda(q)$ of the generic representations of $H_n(q)$
(which themselves may be calculated using the methods and formulae of
\cite{KWy90,KWy92,Rm91,Vj91}).

\proclaim Theorem 2. If $S^{(\lambda_1,\lambda_2)}$ is reducible then,
using the notation of Theorem 1,
\beq
\tilde\chi^\lambda_p=
\sum_{j=0}^{[\lambda_2/p]} \chi^{(\lambda_1+jp,\lambda_2-jp)}(q)
-\sum_{j=0}^{[\mu_2/p]} \chi^{(\mu_1+jp,\mu_2-jp)}(q),
\eeq
where $[x]$ is the largest integer less than or equal to $x$.
\par\noindent
Of course, this Theorem may be used to give the dimension of
$D^\lambda_p$ in terms of the dimensions $d^\nu$ of the
irreducible representations of $S_n$.
For example,
\beq
\dim D^{(6,5)}_3=d^{(6,5)}+d^{(9,2)}-d^{(7,4)}-d^{(10,1)}
=132+44-165-10=1.
\eeq

\begin{center}
{\bf 5\ \ Explicit \boldmath $D^{(\lambda_1,\lambda_2)}_p$}
\end{center}

When $q$ is a primitive $p$th root of unity, the irreducible
$H_n(q)$-module $D^{(\lambda_1,\lambda_2)}_p$ may be constructed
along lines similar to the construction of the Specht modules.
In this, the column and Garnir relations are retained, and are supplemented
by additional relations. These relations will be described elsewhere.
A basis for $D^{(\lambda_1,\lambda_2)}_p$ is defined in terms
of a certain subset of the set of standard tableaux.
In order to specify this set, let
\beq
T^{(\lambda_1,\lambda_2)}=
\youngt{a_1&a_2&a_3&\cdot&\cdot&\cdot&\cdot&\cdot&a_{\lambda_1}\cr
b_1&b_2&b_3&\cdot&\cdot&b_{\lambda_2}\cr},
\eeq
and say that $T^\lambda$ is $s$-strip standard at the $i$th position
if:
\beq
b_i<a_{i+s-2}
\eeq
(or if $i>\lambda_1-s+2$, when of course $a_{i+s-2}$ is undefined).

\proclaim Definition 1. If $\lambda=(\lambda_1,\lambda_2)$ and the
positive integers $p$ and $k$ are such that
$\lambda_1-\lambda_2+2\le kp\le\min\{\lambda_1+1,\lambda_1-\lambda_2+p\}$,
then $T^\lambda$ is said to be {\em $p$-root standard} if $T^\lambda$
is standard and either:
\parindent=30pt
\titem{1. } $T^\lambda$ is $kp$-strip standard at positions
$1,2,\ldots,\lambda_2$;
\titem{or 2. } to the right of the rightmost position of
a non-standard $kp$-strip, there is a position at which
$T^\lambda$ is $((k-1)p+2)$-strip standard.

\noindent
Note that in the important case of $k=1$, the second condition here
can never be satisfied because standardness denies 2-strip standardness.
In this case, the tableaux are identical to those defined in
\cite{Wn88} for the corresponding representations.

As an example, consider $\lambda=(7,4)$ and $p=3$ (so that $k=2$).
In this case, the following are $3$-root standard:
\beq
\youngt{1&2&3&4&6&8&10\cr 5&7&9&11\cr},
\qquad
\youngt{1&3&4&5&6&7&11\cr 2&8&9&10\cr},
\qquad
\youngt{1&2&3&4&5&9&10\cr 6&7&8&11\cr};
\qquad
\label{std2}
\eeq
whereas the following are not $3$-root standard:
\beq
\youngt{1&3&5&6&7&8&9\cr 2&4&10&11\cr},
\qquad
\youngt{1&3&4&5&6&7&10\cr 2&8&9&11\cr},
\qquad
\youngt{1&2&3&4&5&8&10\cr 6&7&9&11\cr}.
\qquad
\label{std3}
\eeq

\proclaim Theorem 3. If $\lambda=(\lambda_1,\lambda_2)$,
the dimension of $D^\lambda_p$ is equal to the number of
$p$-root standard tableaux of shape $\lambda$.

\noindent
The non-generic $H_n(q)$-module $D^\lambda_p$ may then be explicitly
constructed with basis:
\beq
\{v_{t^\lambda}: t^\lambda {\rm\ is\ } p {\rm-root\ standard}\}.
\eeq

\end{document}